\documentclass[reqno, 11pt]{amsart}

\usepackage{amsmath,amssymb,amsthm,amsfonts, amscd, url}
\usepackage{appendix}
\usepackage[numeric]{amsrefs}
\usepackage{hyperref}

\usepackage[dvipsnames]{xcolor}
\usepackage{graphics}
\usepackage{graphicx}

\usepackage{graphicx}
\usepackage{yhmath}
\usepackage{mathdots}
\usepackage{MnSymbol}

\usepackage{ytableau, mathdots}

\newtheorem*{nthm}{Theorem}

\newtheorem*{nlem}{Lemma}

\newtheorem*{nrem}{Remark} 

\newcommand{\ZZ} {\mathbb{Z}}
\newcommand{\CC} {\mathbb{C}}
\newcommand{\RR}{\mathbb{R}}
\newcommand{\QQ}{\mathbb{Q}}

\newcommand{\Oin}{\mathcal{O}}

\newcommand{\la}{\lambda}
\newcommand{\La}{\Lambda}
\newcommand{\Cr}{{\mathcal{C}}}

\newcommand{\g}{{\mathfrak{g}}}

\newcommand{\wt}{{\mathrm{wt}}}

\newcommand{\w}{{\underline{w}}}
\newcommand{\SO}{{\mathrm{SO}}}

\newcommand{\Li}{{\mathfrak{L}}}
\newcommand{\s}{{\mathfrak{s}}}
\newcommand{\x}{{\mathbf{x}}}
\newcommand{\n}{{\mathfrak{n}}}

\newcommand{\cha}{{\mathrm{char}}}
\newcommand{\SL}{{\mathrm{SL}}}
\newcommand{\Sp}{{\mathrm{Sp}}}

\newcommand{\PGL}{{\mathrm{PGL}}}

\newcommand{\W}{{\mathcal{W}}}

\newcommand{\G}{{\mathbf{G}}}

\title{Crystal constructions in Number Theory}
\author{Anna  Pusk\'{a}s}
\date{\today}
\address{University of Massachusetts, Amherst \\ Department of Mathematics and Statistics\\ Amherst, Massachusetts}
\email{puskas@math.umass.edu}

\begin{document}

\maketitle

\begin{abstract}
Weyl group multiple Dirichlet series and metaplectic Whittaker functions can be described in terms of crystal graphs. We present crystals as parameterized by Littelmann patterns and we give a survey of purely combinatorial constructions of prime power coefficients of Weyl group multiple Dirichlet series and metaplectic Whittaker functions using the language of crystal graphs. We explore how the branching structure of crystals manifests in these constructions, and how it allows access to some intricate objects in number theory and related open questions using tools of algebraic combinatorics. 
\end{abstract}

\section{Introduction}\label{sect:introduction}

Crystal graphs are combinatorial objects appearing in the representation theory of semisimple Lie algebras. To an irreducible representation of a semisimple Lie algebra $\g$ one may associate a crystal graph $\Cr .$ The vertices of this graph are in bijection with a weight basis of the representation, and the edges are colored by a set of simple roots of $\g .$ 

Crystals were first studied in connection with the representation theory of the quantized universal enveloping algebra. However, in this chapter it is their structure as a colored (directed) graph and their symmetries related to the Weyl group of $\g$ that are of interest to us. Crystals turn out to be a valuable tool in constructing certain objects from number theory: coefficients of multiple Dirichlet series and metaplectic Whittaker functions. 

Interest in multiple Dirichlet series and metaplectic Whittaker functions is motivated by hard questions in analytic number theory, for example the Lindel\"{o}f Hypothesis, and the study of automorphic forms \cite{chinta2006multiple}. The relevant literature in number theory is extensive (see \ref{subsect:review_lit}). However, since these objects have constructions that are almost purely combinatorial in nature, their study can be approached using tools of algebraic combinatorics.

In this chapter we are interested in this approach. Our main goal is to present combinatorial constructions of metaplectic Whittaker functions and coefficients of multiple Dirichlet series corresponding to root systems of the four infinite families of Cartan types. To do so, we use the language of Littelmann patterns. We highlight how the branching structure of crystals is apparent in the constructions, and indicate how this aspect turns out to be significant in the study of the related questions from number theory. 

Before giving an overview of the structure of the chapter, we say a few more words on the relevant objects. 

Crystal graphs can be parameterized (without referring to the representation theory of the quantum group) using a variety of combinatorial devices, such as the Littelmann path model, Gelfand-Tsetlin patterns, Lusztig's parametrization \cite{lusztig1990canonical, lusztig1991canonical}, tableaux of Lakshmibai and Sheshadri \cite{Lakshmibai-Sheshadri} or Kashiwara and Nakashima \cite{Kashiwara-Nakashima}. For a thorough introduction to the theory of crystals from a combinatorial perspective, the reader is encouraged to consult \cite{bump2017crystal}. 

Here we present crystals in terms of Berenstein-Zelevinsky-Littelmann paths and Littelmann patterns \cite{littelmann1998cones}. Our reason for this choice is twofold. First, most of the constructions in number theory that we are concerned with were either originally given in this language, or are easily rephrased in such terms. Moreover, phrasing the constructions in terms of Littelmann patterns highlights the role of the branching structure of crystals (as well as the significance of some ``nice elements'' of the Weyl group) very well. 

A major hurdle any expository writing on this topic has to overcome is the inherent intricacy and volume of the theory of multiple Dirichlet series and metaplectic Whittaker functions. Since we wish to take a purely combinatorial approach, we largely try and circumvent this issue. {Some background on multiple Dirichlet series and Whittaker functions (as well as on metaplectic groups) will be given in \ref{subsect:MDS_WF}.} For now we say that through their connection to an algebraic group over a local or global field, these objects from number theory are related to the representation theory of the underlying Lie algebra $\g$. Their constructions involve producing, for a dominant weight $\lambda ,$ a polynomial $P_{\lambda} (\x)$ in $r$ variables, where $r$ is the rank of the Lie algebra $\g .$ In section \ref{subsect:MDS_WF} we shall say more about how a Weyl group multiple Dirichlet series or a metaplectic Whittaker function gives rise to such a polynomial $P_{\lambda }.$ However, for most of the chapter we shall ignore details of this background, and concern ourselves with producing a polynomial $P_{\lambda }$ as a sum over a crystal graph. {The ``constructions'' mentioned throughout the chapter refer to constructions of a polynomial $P_{\lambda }$ - depending on the context, this may agree with a $p$-part of a Weyl group multiple Dirichlet series or certain values of a spherical Whittaker function.}

{The combinatorial perspective of focusing our attention on the polynomials $P_{\lambda }$ is, on the one hand helpful when considering questions motivated by the analytic background. On the other hand, these objects are interesting in their own right. This is due to the fact that they can be thought of as deformations of highest weight characters. As a result, techniques of character theory come into play. As a method to study these polynomials it provides insight into the original analytical objects. Furthermore, it motivates further questions.}

{To provide an example we mention two aspects of the polynomials $P_{\lambda }$ now.} One is Weyl symmetry: the polynomials $P_{\lambda }$ inherit certain functional equations under the Weyl group corresponding to the underlying root system. Hence one may construct such a polynomial by (1) taking a sum over an object that is similarly symmetric, such as a highest weight crystal, or (2) by taking an ``average'' over the Weyl group. (See Section \ref{subsect:review_lit} for relevant results in the literature.) Understanding the relationship between these two approaches is a large part of the motivation between studying the constructions combinatorially, and in the cases where the question is resolved, the branching structure of crystals turns out to play a significant role. 

{We briefly explain the relevance of character theory. In the simplest special case $P_{\lambda }$ looks very similar to the deformation of a Schur polynomial; more generally, to the deformation of a Weyl character. (For $P_{\lambda}$ a Whittaker function, this is a consequence of the Casselman-Shalika formula.) Hence one expects that the behavior of families of polynomials $P_{\lambda}$ will be similar to the behavior of characters. On the one hand, this means that identities of Weyl characters may provide a useful tool of study. These come in a couple of different flavors. For example the Weyl(-Kac) character formula produces a character. Branching rules describe the behavior of characters under restriction. Indeed, Tokuyama's theorem (a deformation of the Weyl character formula) turns out to be key in investigating the relationship of the two approaches (1) and (2) mentioned above. Generalizing it to the polynomials $P_{\lambda }$ requires understanding the branching properties of the $P_{\lambda }.$ We shall elaborate on this point in Sections \ref{subsect:Tokuyama_DL} and \ref{subsect:motivation_branching} below. 

On the other hand one may ask if products of the $P_{\lambda }$ satisfy some ``enhanced'' version of other character identities. For example, does a product of such polynomials satisfy ``deformed'' Pieri and Littelwood-Richardson rules? Question of this flavor may be investigated using any description of these objects. We shall see that the $P_{\lambda }$ can be defined in terms of the combinatorial structure of a crystal and a few Gauss sums. Hence any question about them can be phrased in terms of the crystal structure and identities of Gauss sums.
 }

We give an overview of the structure of the chapter. 

\subsection{Structure of the paper}
In the remainder of this Introduction, we first give a brief review of results constructing Weyl group multiple Dirichlet series or metaplectic Whittaker functions (section \ref{subsect:review_lit}). In Section \ref{subsect:Tokuyama_DL} we explain how a theorem of Tokuyama is related to this topic, and how Demazure-Lusztig operators can be used to study, and extend the constructions discussed in this chapter to greater generality. Section \ref{subsect:motivation_branching} provides some further insight into the meaning and significance of branching. 

Littelmann patterns and their bijection with crystal elements are discussed in Section \ref{sect:Littelmann_patterns}. The constructions of Whittaker functions and prime power coefficients of multiple Dirichlet series in terms of highest weight crystals are presented in section \ref{sect:constructions}, and the relationship of the constructions with the branching structure of crystals is highlighted in Section \ref{sect:branching_properties}. 

Section \ref{sect:preliminaries} serves to present some preliminaries. We introduce notation (\ref{subsect:notation}) and present Gauss sums, a necessary arithmetic ingredient to the constructions (\ref{subsect:Gauss_sums}). We then give a brief introduction to crystals and Berenstein-Zelevinsky Littelmann paths (\ref{subsect:hwcr_cones}). We also provide a little more insight into how coefficients of multiple Dirichlet series and Weyl group metaplectic Whittaker functions give rise to polynomials $P_{\la }(\x),$ related to sums over highest weight crystals (\ref{subsect:MDS_WF}).

\subsection{A review of literature}\label{subsect:review_lit}

We discuss the literature of constructions of multiple Dirichlet series and Whittaker functions. Our interest here is from the perspective of combinatorics. Hence we shall focus on the role of the branching structure of the crystals and the significance of special words in the Weyl group. For an insightful and thorough introduction to the topic from a number theoretic perspective, the reader is encouraged to consult \cite{bump2012introduction}, the Introduction of the volume where many of the constructions discussed in this chapter were published. The role of this section is to provide this topic with a wider context; strictly speaking it is not necessary for the understanding of any of the later parts. 

Brubaker, Bump, Chinta, Friedberg, and Hoffstein \cite{bbcfh-wmd1} introduced Weyl group multiple Dirichlet series (WMDS), series in several complex variables with functional equations governed by a finite Weyl group, corresponding to a root system $\Phi $ of finite type. As mentioned above, there are two separate approaches to how to associate a WMDS to a root system: by taking a sum over a crystal, or by Chinta and Gunnells \cite{cg-jams}, by averaging over the Weyl group. 

The authors of \cite{bbcfh-wmd1} conjecturally related WMDS to Whittaker coefficients of metaplectic Eisenstein series. This connection is of interest in that it allows one to prove functional equations and analytic continuation of the constructed series. Elucidating this connection motivates study of these objects as well. In the following paragraphs when we refer to a ``conjectural description'' of a WMDS as a sum over a crystal, we mean either that the constructed series is conjectured to be the Whittaker coefficient of a metaplectic Eisenstein series, or that it is conjectured to agree with a series constructed via the averaging method. 

We shall mention relevant results in all four infinite families of Cartan types; some of these results will be covered in more detail in section \ref{sect:constructions}. 

Brubaker, Bump and Friedberg \cite{bbf-annals} describe the Fourier-Whittaker coefficients of Eisenstein series on a metaplectic cover of the general linear group as a Weyl group multiple Dirichlet series. They compute the prime power coefficients ($p$-parts) of these series in terms of the string parametrization of a crystal by Berenstein and Zelevinsky \cite{berenstein1996canonical, berenstein_zelevinsky_Tensorproductmult} and Littelmann \cite{littelmann1998cones}. In \cite{bbf-wmdbook} the same authors further explore the combinatorics in the type $A$ case. They give two separate constructions of the $p$-part in type $A$. These can be seen as corresponding to two different choices of {\em{nice decompositions}} of the long element of the Weyl group. The authors then prove that the two descriptions give the same $p$-parts through a subtle combinatorial argument. The equivalence of the two statements allows them to prove analytic continuation and functional equations for the emerging multiple Dirichlet series. In proving the equivalence, they observe the significance of some purely combinatorial phenomena - such as the Sch\"{u}tzenberger involution. Their method provides an example of how to build $p$-parts of multiple Dirichlet series out of finite crystal data. 

Beineke, Brubaker and Frechette \cite{beinekebrubakerfrechette, beineke2012crystal} give a definition for a WMDS in terms of statistics on a highest weight crystal of Cartan type $C$. They prove analytic continuation and functional equation of such series using a connection to Eisenstein series over odd orthogonal groups in the nonmetaplectic case, and conjecture a similar connection in general. 

Friedberg and Zhang \cite{friedberg2015eisenstein} study Eisenstein series over metaplectic covers of odd orthogonal groups. They then describe the $p$-parts of the MDS that are the Whittaker coefficients of these series in terms of type $C$ highest weight crystals. They in fact give two descriptions. The first one is only valid in the case of odd covers; this proves the conjectured connection in \cite{beineke2012crystal} above. The second is uniform in the degree $n$ of the metaplectic cover, but the assignment of number theoretic data to the combinatorial structure is more subtle. An interesting feature of their methods is that they are is inductive by rank. Furthermore the proof of the agreement of two descriptions relies on the type $A$ theory of \cite{bbf-wmdbook}. 

As for type $D,$ Chinta and Gunnells \cite{chinta2012littelmann} give a conjectural construction of a Weyl group multiple Dirichlet series of type $D.$ The $p$-part of a series is produced as a sum over a highest weight crystal associated to an irreducible representation of $SO(2r).$ The contribution of a crystal element to the sum is described in terms of the corresponding Littelmann pattern. 

We also mention constructions of Whittaker functions. McNamara \cite{mcnamara} considers Whittaker functions on metaplectic covers of a simple algebraic group over a nonarchimedean local field. The Whittaker function is given as an integral (over the unipotent radical). Given a reduced decomposition of the long element in the Weyl group, one may break up the domain of integration into a set of cells. These cells are in a natural bijection with elements of an (infinite) crystal. By computing the integral on each cell, the Whittaker function is produced as a sum over a(n infinite) crystal structure. In type $A,$ the resulting formula for the Whittaker function agrees with the formula given for local parts of a Weyl group multiple Dirichlet series by \cite{bbf-annals}. 

The averaging approach can also be used to construct metaplectic Whittaker functions as shown by Chinta and Offen \cite{co-cs} in the type $A$ case, and McNamara \cite{mcnamara2} in general. Work of McNamara thereby provides a number-theoretic proof that the two methods (averaging and crystal constructions) produce the same local parts. 

From a combinatorial perspective, the formulas produced by the two separate approaches (averaging, or sum over a crystal) are related in the nonmetaplectic case by a theorem of Tokuyama \cite{tokuyama-generating}. More generally, Demazure operators can be used to elucidate the connection between the two approaches \cite{cgp, puskas2016whittaker} combinatorially. This relies heavily on the branching properties of crystals. This is explained in more detail in sections \ref{subsect:Tokuyama_DL} and \ref{subsect:motivation_branching}. 

{The reason for the emergence of crystal bases in the study of these topics in itself warrants further exploration. Some results of this flavor exist both in the local and in the global setting. As mentioned above, the work of McNamara \cite{mcnamara} expresses a Whittaker function as a sum over cells of the unipotent radical. The cell decomposition is then related to geometric realizations of the crystal in terms of Lusztig data \cite{lusztig1996algebraic} and Mirkovi\'c-Vilonen cycles \cite{braverman2001crystals}. In the global setting Brubaker and Friedberg \cite{brubaker2015whittaker} study Whittaker coefficients of metaplectic Eisenstein series induced from maximal parabolics. They produce a formula for the Whittaker coefficient for a wide class of long words in the Weyl group by matching contributions with Lusztig data through MV polytopes considered by Kamnitzer \cite{kamnitzer2010mirkovic}.}

In addition, highest weight crystals are not the only combinatorial device that is of use in constructing these number theoretic objects. An other approach uses metaplectic ice \cite{2010arXiv1009.1741B, brubaker2016schulzarxiv, brubaker2016bumpbuciumas}.

\subsection{Tokuyama's Theorem}\label{subsect:Tokuyama_DL}

In this section we explain how the results above relate to a deformation of the Weyl character formula by Tokuyama \cite{tokuyama-generating}, and how understanding the branching structure of crystals elucidates the relationship of the constructions. 

The constructions produce a polynomial $P_{\la }(\x)$ that satisfies certain functional equations under a Weyl group $W.$ (Here we assume that the polynomial ring $\CC[\x]$ in $r$ variables is identified with the group algebra of the weight lattice.) This is done either by taking a sum over a highest weight crystal $\Cr _{\la }$, or by taking an average over the Weyl group. We explain how both of these strategies results in a polynomial that is, roughly speaking, a deformation of a Weyl character. 

First let us consider the method of producing the polynomial $P_{\lambda }$ by taking a sum over the crystal graph:
\begin{equation}\label{eq:crystal_sum_roughshape}
P_\la(\x)=\sum_{b\in \Cr_{\la}}G(b)\x^{\wt(b)}
\end{equation}
The elements of the crystal $\Cr_{\la}$ are in bijection (via the weight map $\wt$) with a weight basis of a representation of highest weight $\la $. Note that if we had $G(b)=1$ for every element of the crystal in \eqref{eq:crystal_sum_roughshape}, then the resulting sum would be the character $\chi_{\la }$ of this highest weight representation:
\begin{equation}\label{eq:char_twoforms}
\sum_{b\in \Cr_{\la}} 1\cdot \x^{\wt(b)} = \chi _{\la }(\x)=\frac{1}{\Delta }\cdot \sum_{w\in W} (-1)^{\ell(w)} \cdot w(\x^{\la +\rho })
\end{equation}
Here the right-hand side is the Weyl character formula, and $\Delta $ is a Weyl denominator. 

In general $G(b)$ is more complicated, but in the simplest case we have that in fact 
\begin{equation}\label{eq:tokuyama_roughshape}
P_\la(\x)=\sum_{b\in \Cr_{\la}}G(b)\x^{\wt(b)}=\Delta_q\cdot \chi _\la(\x)
\end{equation}
where $\Delta _{q}$ is a deformation of the Weyl denominator.

Next let us consider the ``averaging approach'' to constructing the polynomial $P_{\la }(\x).$ This approach produces $P_{\la }(\x)$ by an expression similar to the right-hand side of \eqref{eq:char_twoforms}. However, the action of $w$ on the monomial $\x^{\la +\rho }$ is replaced by the Chinta-Gunnells action \cite{cg-jams}, which depends on the metaplectic degree $n.$ In the special case of $n=1,$ this construction results in the expression $\Delta_q\cdot \chi _\la(\x)$ as well. 

The statement that the two approaches to constructing $P_{\lambda }$ give the same result can thus be phrased as a combinatorial identity, a deformation of the Weyl character formula.
In the nonmetaplectic case, this is the second equality in \eqref{eq:tokuyama_roughshape}, and this identity is a theorem by Tokuyama \cite{tokuyama-generating}. 

When $n>1,$ then understanding the relationship between the two constructions of $P_\la $ combinatorially amounts to proving a metaplectic analogue of Tokuyama's theorem. In the type $A$ case, this was done by the author in \cite{puskas2016whittaker} using metaplectic Demazure-Lusztig operators defined in \cite{cgp}. 

We mention that analogues of Tokuyama's theorem for root systems of other types have been given by Hamel and King \cite{Hamel-King} (for type $B$) and Friedlander, Gaudet and Gunnells \cite{Friedlander-Gaudet-Gunnells-typeG2} (in type $G_2$). Note also that the agreement between the relevant constructions in the type $A$ case follows from work of McNamara as indicated above. However, treating these sides combinatorially via Demazure-Lusztig operators allows one to understand how the constructions can be extended to more general settings, for example, from the finite dimensional to the general Kac-Moody setting \cite{patnaik2017metaplectic}, or, from Whittaker functions to the constructions of Iwahori-Whittaker functions \cite{patnaik2017iwahori}.

The proof of the metaplectic analogue of Tokuyama's theorem in \cite{puskas2016whittaker} relies heavily on the type $A$ crystal construction ``respecting'' the branching structure of the highest weight crystal. We explain this in more detail next. 

\subsection{Motivation: Demazure-Lusztig operators and the branching structure}\label{subsect:motivation_branching}

As seen above, understanding the combinatorial relationship between different constructions of the polynomial $P_{\la }$ (which may be a Whittaker function or the prime part of a WMDS) amounts to proving a metaplectic analogue of Tokuyama's theorem. 

Using Demazure-Lusztig operators one may phrase a more general identity, corresponding to elements $w$ of the Weyl group, and any metaplectic degree $n$. The more general identity \cite[Theorem 1.]{puskas2016whittaker} is of the form:
\begin{equation}\label{eq:gentok_roughshape}
\left(\sum_{u\leq w} T_u\right) \x ^{\la } = \sum _{b\in \Cr _{\la }^{(w)}} G(b) \x ^{\wt (b)}
\end{equation}
Here the expression on the left-hand side can be thought of as the general form of the expression produced by the averaging method (by results in \cite{cgp}) and the right-hand side is a sum over a Demazure crystal. The ``metaplectic analogue of Tokuyama's theorem'' (in type $A$) is the special case of this statement corresponding to $w$ being the long element of the Weyl group. This more general statement has the advantage that it can be proven ``one simple reflection at a time,'' i.e. by induction on the length $\ell(w)$ of the Weyl group element. 

The fact that the construction of $P_{\la }$ as a sum over a highest weight crystal respects the branching structure of the crystal is crucial to the proof. We explain what we mean by this below. We shall return to this discussion in more detail in section \ref{sect:branching_properties} equipped with the necessary background. 

The crystal $\Cr_{\la }$ is graph, whose edges are labeled by simple roots $\alpha _i$ ($1\leq i\leq r$) of an underlying Lie algebra or rank $r.$ When the edges labeled by $\alpha _r$ are omitted, the remaining graph is a disjoint union of crystals $\Cr _{\mu }$, corresponding to a Lie algebra of the same Cartan type as $\g ,$ but rank $r-1:$ 
\begin{equation}\label{eq:crystal_cranching_intro}
\Cr _{\la }=\sqcup _{\mu }\Cr _{\mu }
\end{equation}
A crystal element $b\in \Cr_{\la }$ has a contribution $G(b)=G_{\la }(b)$ in the sum \eqref{eq:crystal_sum_roughshape}. The element $b$ belongs to exactly one of the rank $r-1$ crystals $\Cr_{\mu }\subset \Cr _{\la }.$ The element $b$ has a contribution in the analogous construction of $P_{\mu }(\x).$ By the constructions respecting the branching structure, we mean that we have 
\begin{equation}\label{eq:contrib_decompose}
G_\la (b)=g(\mu) \cdot G_\mu (b),
\end{equation} 
where the factor $g(\mu) $ is the same for every element $b\in \Cr_{\mu}\subset \Cr _{\la }.$ It follows that $P_{\la }$ can be written as an expression of polynomials $P_{\mu }$ corresponding to the weights $\mu $ of the decomposition \eqref{eq:crystal_cranching_intro} above. 

This means that statements about these crystal constructions are amenable to proof by induction on rank. The parametrization of crystal elements by Littelmann patterns highlights the branching structure of crystals. We encourage the reader to keep the branching structure in mind while reading through the sections covering the ingredients of the constructions. 

\subsection{Acknowledgments}
I would like to thank the editors of this volume for giving me the opportunity to contribute. I am grateful to several people for helpful conversations and advice during the writing of this chapter, including Holley Friedlander, Paul E Gunnells, Dinakar Muthiah and Manish Patnaik. During parts of the writing process, I was a postdoctoral fellow at the University of Alberta and a visiting assistant professor at the University of Massachusetts, Amherst, and I am grateful to both institutions. While at the University of Alberta, I was supported through Manish Patnaik's Subbarao Professorship in number theory and an NSERC Discovery Grant. {I also thank the referees for their insightful comments for the improvement of this chapter. In particular, we thank one of the referees for their comments on the connections to character theory and on Gauss sums.}

\section{Preliminaries}\label{sect:preliminaries}

Before describing the constructions mentioned above, we cover a few preliminaries. The constructions in section \ref{sect:constructions} have two main ingredients: a set of root data, and an arithmetic ingredient in the form of certain Gauss sums. We introduce notation and the necessary ingredients below. 

\subsection{Notation}\label{subsect:notation}

Throughout the paper, $\Phi $ shall denote a root system of rank $r$, with $\Phi^+$ (respectively $\Phi^-$) being the set of positive (respectively, negative) roots. Let $\Delta=\{\alpha_1, \ldots ,\alpha _r\}$ be a set of simple roots in $\Phi ,$ let us write $\rho =\frac{1}{2}\sum_{\beta \in \Phi ^+}\beta $ for the Weyl vector. Of particular interest are the root systems of Cartan types $A,$ $B,$ $C$ and $D.$ We give \cite{bourbaki} as a general reference on root systems. Note however that when discussing Littelmann patterns \cite{littelmann1998cones}, our numbering of the simple roots agrees with that of Littelmann and hence differs from that of Bourbaki. (For example in type $D$ the simple roots $\alpha _1$ and $\alpha _2$ are orthogonal.)

Let us write $\sigma_1,\ldots ,\sigma _r$ for the set of simple reflections corresponding to the simple roots; $\sigma_i$ is the reflection through the hyperplane perpendicular to $\alpha_i.$ The Weyl group $W$ is generated by the simple reflections $\sigma_i$ ($1\leq i \leq r$). Every element $w\in W$ can be written as a product $w=\sigma_{i_1}\cdots \sigma_{i_k}.$ We call this a reduced decomposition and $\w = [i_1,\ldots ,i_k]$ a reduced word if $k$ is minimal and $k=\ell(w),$ the length of $w.$ The Weyl group has a unique longest element $w_0\in W.$ The parameterization of highest weight crystals by Littelmann patterns given in section \ref{sect:Littelmann_patterns} depends on a choice of a {\em{nice decomposition}} $\w_0$ of the long element. 

The Weyl group permutes the elements of $\Phi .$ Let $\Phi(w)= w^{-1}(\Phi^-)\cap \Phi ^+ ,$ then $\ell(w)=|\Phi(w)|$ and $\Phi(w_0)=\Phi^+.$ We shall denote the weight lattice corresponding to $\Phi $ by $\Lambda ,$ and the fundamental weights corresponding to the basis $\Delta $ by $\varpi _1, \ldots ,\varpi _r.$ The constructions we are concerned with produce polynomials in $\CC[\Lambda ].$ The Weyl group has a natural action on $\La $ and hence on $\CC[\La ],$ it thus make sense to talk about the resulting polynomials being symmetric (or having functional equations) under the Weyl group. 

\subsection{Gauss sums}\label{subsect:Gauss_sums}

Next we introduce notation for the arithmetic ingredients of the polynomials constructed in Section \ref{sect:constructions}. The contribution of a crystal element is given by $n$th order Gauss sums, where $n$ is a positive integer. (In applications, $n$ is the degree of the metaplectic cover.) Recall that (for $n=2$) one may take the quadratic Gauss sum 
\begin{equation}\label{eq:quad_Gauss}
G(a)=\sum_{k=0}^{p} \left(\frac{k}{p}\right)e^{\frac{ka2\pi i}{p}}
\end{equation}
where $\left(\frac{k}{p}\right)$ is the Legendre symbol, i.e it is $1$ if $k$ is a square modulo $p$ and $-1$ otherwise. The Gauss sums appearing in the construction are generalizations of the one in \eqref{eq:quad_Gauss}. The Legendre symbol is replaced by an $n$th power residue symbol (a multiplicative character), and $e^{\frac{ka2\pi i}{p}}$ is replaced by an additive character. 

To state the definition of the Gauss sums $g(a)(=g_t(a))$ and $h(a)(=h_t(a))$ we introduce some more notation. We use the language of local and global fields, and provide examples (see \cite{neukirch2013algebraic} for reference). The Gauss sums $g(a)$ and $h(a)$ are functions that depend on the residue of $a$ modulo $n.$ For the purpose of understanding the constructions and their relationship to the branching structure of crystals, the values of these functions is not crucial.

Following \cite{bbf-wmd2} let $F$ be a global field; the reader may choose to simply think of $\QQ$ as an example. For a place $v$ of $F$ one may take the completion $F_v.$ (For example, the completion $\QQ_p$ of $p$-adic numbers at any finite prime $p$, or the completion $\RR$ at the infinite place.) Let $\Oin_v$ denote the set of integers (e.g. $\ZZ_p\subset \QQ_p$ or $\ZZ\subset \RR$). Let $S$ be a finite set of places of $F,$ and let $\Oin_S$ denote the set of $S$-integers $x\in K$ such that $x\in \Oin _v$ for every $v\notin S.$ (For $S=\{\infty,2\}$ the set $\Oin_S\subseteq \QQ$ is the set of rational numbers with only $2$ in the denominator.) For a sufficiently large $S,$ $\Oin_S$ is a principal ideal domain. Let $F_S=\prod_{v\in S}F_v,$ $\Oin_S$ embeds into $F_S$ diagonally. Let $\psi $ be a character of $F_S$ trivial on $\Oin _S$ and no larger fractional ideal. Let $\left(\frac{\cdot }{\cdot }\right)_n$ denote the $n$-th order residue symbol and $t$ a positive integer. We define 
\begin{equation}\label{eq:def:gt}
g_t(a,c)=\sum_{\substack{d \mod c}} \left(\frac{d}{c}\right)_n^t\psi \left(\frac{ad}{c}\right) .
\end{equation}
The constructions in section \ref{sect:constructions} will involve special values of $g_t(a,c).$ We shall have $t=1$ or $t=2$ be the length of a simple root (i.e. $t=1$ in the simply laced cases and $t=1$ or $t=2$ in type $B$ or type $C$), and we shall have fixed a prime $p.$ Then we set 
\begin{equation}\label{eq:gt_ht}
g_t(a)=g_t(p^{a-1},p^a);\text{ and }h_t(a)=g_t(p^a,p^a)=\left\lbrace \begin{array}{ll}
|(\Oin _S/p\Oin_S)^{\times }| & \text{if }t^{-1}n\mid a\\
0 & \text{if }t^{-1}n\nmid a
\end{array}\right.
\end{equation}

In the remainder of the paper we use the notation $q=|\Oin_S/p\Oin _S|$ for the order of a residue field. 

{The polynomials $P_{\lambda }$ that we shall define in Section \ref{sect:constructions} are given as a sum over a crystal. Each term is determined via combinatorial data coming from the parameterization of the corresponding crystal element via a method that makes use of the above Gauss sums. For example in Cartan type $A$ the Gauss sums take $n$ values (indexed by the residue classes modulo $n$). Consequently any statement about these polynomials can be phrased entirely in terms of the structure of the highest weight crystal and identities of these $n$th order Gauss sums. Such identities are rare. In addition to the identity expressing the relationship of Gauss sums corresponding to conjugate characters, one has the Hasse-Davenport relations \cite{davenport1934nullstellen}. By work of Yamamoto \cite{yamamoto1966conjecture} these are essentially the only multiplicative identities of these Gauss sums. We thank one of the referees for pointing this out.}

\subsection{Highest weight crystals and Littelmann's cone}\label{subsect:hwcr_cones}


Given an irreducible (finite) root system $\Phi $ and a dominant weight $\la $ there is an associated crystal graph $\Cr_{\la }$. We shall describe the structure of $\Cr_{\la }$ as a directed graph with colored edges. 
We mention that if $\g$ is a simple Lie algebra with root system $\Phi $, and $V_{\la }$ is the unique simple $\g$ module with highest weight $\lambda $ then the quantized universal enveloping algebra $U_{\mathsf{q}}(\g)$ has a corresponding module. A crystal base is a base for this module at ${\mathsf{q}}=0.$ It carries a graph structure induced by its structure as a $U_{\mathsf{q}}(\g)$ module. For further information, see \cite{kashiwara1995crystal}. 
Here we forgo exploring the connection with the quantum group. We instead explain the structure of a crystal as a colored directed graph and the parameterization of crystals by Berenstein-Zelevinsky-Littelmann paths and Littelmann patterns. 

\subsubsection{The crystal as a colored directed graph}

We now describe $\Cr_{\la}$ as a colored directed graph. Let $B$ be a finite set, we call elements of $B$ elements (or vertices) of the crystal. (We shall abuse notation and write $b\in \Cr_{\la }$ for a $b\in B$.) For every $1\leq i\leq r$ we have operators $f_i:B\cup \{0\}\rightarrow B\cup \{0\}$ and $e_i:B\sqcup \{0\}\rightarrow B\sqcup \{0\}$ acting on the vertices. We shall refer to these as root operators. They have the property that if $b,b'\in B$ then $f_ib=b'$ and $b=e_ib'$ are equivalent. This defines the structure of $\Cr_\la $ as a colored directed graph: if $b,b'\in B$ and $f_ib=b',$ then $\Cr _{\la }$ has a directed edge $b\stackrel{i}{\rightarrow}b'$ ``colored'' by the index $i.$ There is a weight function $\wt:B\rightarrow \Lambda $ such that $\wt(f_i(b))=\wt(b)-\alpha _i$ and in fact the function $\wt $ is a bijection between $B$ and a weight basis of the highest weight $\g$-module $V_{\lambda }.$ In particular, there is a unique ``highest element'' $b_{\lambda }\in \Cr _{\lambda }$ with $\wt(b_{\lambda })=\lambda .$ This $b_{\la}$ is the unique element of $B$ such that $e_ib_{\la }=0$ for every $1\leq i\leq r.$ It follows that 
\begin{equation}\label{eq:crystal_is_reachable_from_highest}
B\sqcup\{0\}=\{f_{i_1}^{n_1}f_{i_2}^{n_2}\cdots f_{i_k}^{n_k}b_{\lambda } \mid 1\leq i_j\leq k, 0\leq n_j\}.
\end{equation}
We shall be interested in writing an element $b\in B$ as $b=f_{i_1}^{n_1}f_{i_2}^{n_2}\cdots f_{i_k}^{n_k}b_{\lambda } $ in particular when the sequence of indices $[i_1, i_2\ldots,i_k]$ is a reduced word. 

\subsubsection{Berenstein-Zelevinsky-Littelmann paths}\label{subsubsect:BZLpaths}
Let $\w=[i_1, i_2\ldots,i_k]$ be a reduced word in $W$ and let $\n=[n_1,n_2\ldots ,n_k]\in (\ZZ_{\geq 0})^k$ for $b=f_{i_1}^{n_1}f_{i_2}^{n_1}\cdots f_{i_k}^{n_k}b_{\lambda }.$ We call $\n$ an adaptive string of $b$ \cite{littelmann1998cones} if for every $1\leq j\leq k$ we have 
\begin{equation}\label{eq:adaptive_string}
1\leq j\leq k:\ e_{i_j}f_{i_{j+1}}^{n_1}\cdots f_{i_k}^{n_k}b_{\lambda }=0
\end{equation}
We can think of $\n$ as encoding a path from $b$ to $b_{\la }$ along crystal edges (against the direction of the edges), using $\w$ as a road map. To get the path, starting at $b$ we first take steps along edges colored $i_1$ as long as that is possible. After taking $n_1$ steps, we arrive at a vertex $b_1=e_{i_1}^{n_1}b$ such that $e_{i_1}b_1=0.$ We then proceed with steps along edges colored $i_2$ for as long as possible, etc. 

Taking an adaptive string above defines a map $b\mapsto {\ZZ_{\geq 0}}^{\ell(\w)}$ for any reduced word $\w.$ Let $\w_0$ be a long word of the Weyl group. Write $S_{\w_0}^{\la }\subseteq {\ZZ_{\geq 0}}^{\ell(\w_0)}$ for the set of adaptive strings that occur in $\Cr_{\la }$ and $S_{\w_0}\subseteq {\ZZ_{\geq 0}}^{\ell(\w_0)}$ for the set of strings that occur for any strongly dominant $\la .$ Then it follows from work of Littelmann, Berenstein and Zelevinsky \cite{berenstein1993string, littelmann1994littlewood, littelmann1998cones} that $S_{\w_0}$ is the set of integral points inside a convex cone, which we from now on refer to as the {\em{Littelmann cone}} $C_{\w_0}.$ Furthermore, the set $S_{\w_0}^{\la }$ is the set of integral points in a convex polytope $C_{\w_0}^{\la }$ in this cone (the {\em{Littelmann polytope}}). The inequalities describing $C_{\w_0}$ depend on the long word $\w_0 ;$ the further inequalities describing $C_{\w_0}^{\la }$ depend on $\lambda $ as well. For particularly ``nice'' choices of $\w_0$ \cite{littelmann1998cones} these inequalities take on a transparent form. We shall describe these choices for Cartan types $A,$ $B,$ $C$ and $D$ as well as the Littelmann patterns they give rise to in section \ref{sect:Littelmann_patterns}. For $\w=\w_0$ we shall refer to the adaptive string $\n$ corresponding to a vertex $b\in \Cr_{\la }$ (as well as the corresponding path in $\Cr_{\la }$) as the Berenstein-Zelevinsky-Littelmann path or $BZL$ path of $b$ and write $BZL(b)=\n .$ 

\subsection{Multiple Dirichlet series and Whittaker functions}\label{subsect:MDS_WF}

We briefly introduce the objects from number theory that are produced by the constructions in Section \ref{sect:constructions}. Since we wish to focus on the combinatorics of the constructions, we keep the length of this section to a minimum. Our purpose here is merely to motivate the appearance of highest weight crystals as an apt combinatorial device in the study of these objects. 

\subsubsection{Multiple Dirichlet series}\label{subsubsect:MDS}

Multiple Dirichlet series are series in several complex variables. They can be used to study automorphic $L$-functions, generalizations of the Riemann zeta function via the Langlands-Shahidi method. Of special interest to us here are Weyl group multiple Dirichlet series, whose functional equations are governed by a Weyl group associated to a (finite) Cartan type. The functional equations are of significance in proving meromorphic continuation and functional equations. We explain briefly how prime power coefficients of multiple Dirichlet series are related to sums over a highest weight crystal. We follow the notation of \cite{bump2012introduction} with some simplifications, so as not to occlude the picture. 

Let $F$ now be a global field. We wish to construct a series in $r$ variables $s_1,\ldots ,s_r$: 
\begin{equation}\label{eq:Dirichle_series}
\sum_{C_i} H(C_1,\ldots ,C_r; m_1,\ldots ,m_r) \cdot |C_1|^{-2s_1}\cdots |C_r|^{-2s_r}
\end{equation}
where the summation is over ideals $C_i$ of $\Oin_S.$ Relating such a series to automorphic $L$-functions imposes certain restrictions on its construction. For example, though a series does not have an Euler product in the way the Riemann zeta function does: 
\begin{equation*}\label{eq:Riemann_zeta_Euler_product}
\zeta (s)=\sum_{n=1}^{\infty} \frac{1}{n^s} =\prod_{p \text{ prime}} \frac{1}{1-p^{-s}}
\end{equation*}
its coefficients satisfy a twisted multiplicativity and the series is hence determined by its $p$-parts 
\begin{equation}\label{eq:MDS-ppart}
\sum_{k_i=1}^{\infty } H(p^{k_1},\ldots ,p^{k_r}; p^{l_1},\ldots ,p^{l_r}) \cdot |p|^{-2k_1s_1-\cdots -2k_rs_r}
\end{equation}
where $p$ is (a representative of) a prime ideal, and $(l_1,\ldots ,l_r)$ correspond to a weight $\lambda =\sum_{i=1}^{r}l_i\varpi_i .$ 
 
Constructing a Weyl group multiple Dirichlet series thus amounts to describing the coefficients $H(p^{k_1},\ldots ,p^{k_r}; p^{l_1},\ldots ,p^{l_r})$ for any fixed weight $\lambda.$ Note that assigning a weight to every $(k_1,\ldots ,k_r)$ as above, we may interpret the $p$-part as a sum over the weight lattice $\Lambda .$ Its support turns out to be finite, and in fact contained in the convex hull of the Weyl group orbit of $\lambda .$ 

Recall that for a highest weight crystal $\Cr_{\la }$ (associated to a root system $\Phi $) the weight function $\wt:\Cr_{\la }\rightarrow \La $ is a bijection between vertices of $\Cr_{\la }$ and a weight basis of a representation with highest weight $\lambda .$ Hence the constructions of the $p$-part may be written as a sum over a highest weight crystal.
 
\subsubsection{Whittaker functions}\label{subsubsect:Whittaker_fns}
Our aim here is to motivate why metaplectic analogues of the Casselman-Shalika formula lead to constructions involving highest weight crystals. 

Let $\G$ be a split reductive group defined over $\ZZ .$ (The reader may think of $\SL_{r}$ or $\Sp_{2r}$.) Let $F$ be a nonarchimedean local field (for example $F=\QQ_p,$ the $p$-adic numbers), and $\Oin\subset F$ the ring of integers in $F$ (e.g. $\Oin=\ZZ_p$). Let $G=\G (F)$ and $K=\G(\Oin )$ be a maximal compact in $G.$ Let $T\subseteq G$ be a maximal torus, and $U\subset G$ be the unipotent radical of a Borel subgroup of $G.$ (In the examples above, $T$ is the group of diagonal matrices in $G$ and $U$ the group of upper triangular matrices with $1$s on the diagonal.) Let $\widehat{G}$ denote the Langlands dual of $G$ (we have $\widehat{\SL}_{r+1}=\PGL_r$ and $\widehat{\Sp}_{2r}=\SO_{2r+1}$); let $\Phi $ be the root system associated with $\widehat{G}$ and $\La $ its weight lattice. (Here $\Phi $ is of type $A$ or type $B$ for $\SL_{r}$ or $\Sp_{2r}$ respectively.)
To an element $\x\in \widehat{T},$ we may associate a Whittaker function $\W :G\rightarrow \CC $ that satisfies $\W(ugk)=\psi (u)\W(g)$ (for $u\in U,$ $g\in G$ and $k\in K,$ where $\psi $ is an unramified character of $U$). Let $\pi \in F$ be a uniformizer (e.g. $p$ in $F=\QQ_p$). By the Iwasawa decomposition, any element $g\in G$ can be written as $g=u\pi^{\la }k,$ where $u\in U,$ $k\in K,$ and $\la \in \La $ is a cocharacter of $T.$ A Whittaker function $\W$ is then determined by its values $\W(\pi^{\la }).$ In this classical (nonmetaplectic) setting, these values are determined by the Casselman-Shalika formula \cite{casselman-shalika-cs2}: 
\begin{equation}\label{eq:Casselman-Shalika}
\W(\x,\la)=\prod_{\alpha \in \Phi ^+}(1-q^{-1}\x ^{\alpha }) \chi _{\lambda }(\x)
\end{equation}
This expresses the values of a Whittaker function $\W(\x,\la)$ in terms of the character $\chi _{\lambda }(\x)$ of a representation of $\widehat{G}$ of highest-weight $\la .$ (Here $q=|\Oin/\pi\Oin|$ as before.)

Now let $n$ be a positive integer so that  $\cha F\nmid n$ and $|\mu_{2n}|=2n$ for the group $\mu_{2n}\subset F$ of $n$th roots of unity. Then an $n$-fold metaplectic cover $\tilde{G}$ of $G$ is a central extension 
$$1\rightarrow \mu_n\rightarrow \tilde{G}\rightarrow G\rightarrow 1$$
constructed from the Cartan datum of $G$ and some arithmetic data on $F$ \cite{matsumoto1969sous}. 
By a metaplectic generalization of the Casselman-Shalika formula we mean an analogue of \eqref{eq:Casselman-Shalika} for Whittaker functions on $\tilde{G}$. As mentioned in \ref{subsect:Tokuyama_DL} such a generalization may produce such a formula as a sum over a highest-weight crystal, or as a sum over a Weyl group. This is motivated by the shape of \eqref{eq:Casselman-Shalika}, and its similarity with the deformation of the Weyl character formula in \eqref{eq:tokuyama_roughshape}. 

%

\section{Littelmann patterns}\label{sect:Littelmann_patterns}

We recall Littelmann patterns from \cite{littelmann1998cones} in each of the Cartan types $A_r,$ $B_r,$ $C_r$ and $D_r$. A pattern is an array of $\ell(w_0)$ nonnegative integers. Integral points of the Littelmann cone (see \ref{subsect:hwcr_cones}) are in bijection with the set of patterns that satisfy a set of inequalities. Imposing a further set of inequalities gives a parametrization of integral points within the Littelmann polytope, i.e. a highest weight crystal for a fixed highest weight. The contribution of a single element to the sums in Section \ref{sect:constructions} will be phrased in terms of the corresponding Littelmann pattern. 

The branching properties of highest weight crystals and how it is reflected in the constructions will be made explicit in section \ref{sect:branching_properties}. One may observe these branching properties in the extent to which the Littelmann patterns are consistent within an infinite family of Cartan types. Note also that the simple root $\alpha _r$ that is ``new'' in rank $r$ is associated only to entries in the top row of the pattern. 


\subsection{The choice of a long words}\label{subsect:choice_longword}

Recall that a long word is a reduced decomposition of the long element of the Weyl group. The parametrization of crystal elements in terms of Littelmann patterns is dependent on the choice of a long word $\w_0.$ The choice of particular ``nice'' long words results in the Littelmann cone having a transparent description. We give the nice long words here for each infinite family of Cartan types.

Notice that the choice is consistent within each family in the following sense. Let $X$ stand for any  of $A,$ $B,$ $C$ or $D$ and let $\w_0^{X_r}$ be the choice of long word for type $X_r,$ i.e. rank $r.$ Then the word $\w_0^{X_r}$ starts with the long word $\w_0^{X_{r-1}}$ from rank $r-1.$

The choices are as follows. 

\begin{eqnarray}\label{eq:choice_longword}
\w_0^{A_r} & = & [(1),(2,1),(3,2,1),\ldots ,(r,r-1,\ldots ,2,1)]\\ \label{eq:nice_long_word_A}
\w_0^{B_r} = \w_0^{C_r} & =& [(1),(2,1,2),\ldots ,(r,r-1,\ldots ,2,1,2,\ldots ,r)]\\ \label{eq:nice_long_word_BC}
\w_0^{D_r} & = & [(1),(2),(3,1,2,3),\ldots ,(r,r-1,\ldots ,3,1,2,3,\ldots ,r)] \label{eq:nice_long_word_D}
\end{eqnarray}


\subsection{The shape of patterns}\label{subsect:shape_of_patterns}

The choice of a long word $\w_0$ establishes a bijection between elements of a crystal and $\ell(\w_0)$-tuples of nonnegative integers via BZL paths as in \ref{subsubsect:BZLpaths}. We arrange these $\ell(w_0)$ integers as entries $a_{i,j}$ of a Littelmann pattern. The shape of the pattern reflects the choice of $\w_0$ made. 

Each column of a pattern corresponds to a particular index $1\leq j \leq r.$ Entries $a_{i,j}$ with the same column index $j$ correspond to occurrences of same simple reflection in the word $\w_0$. A row of the pattern will correspond to a step in the rank within the infinite family of Cartan types. 

In the remainder of this chapter we follow the convention that if $a_{i,j}$ is not an entry of a pattern, then $a_{i,j}=0.$ (This is the case for example if $i\leq 0$ or $j<i.$)

%

\subsubsection{Type $A_r$} 
We have $\ell(\w_0^{A_r})-\ell(\w_0^{A_{r-1}})=r$ for $r\geq 2.$ 
A Littelmann pattern of type $A_r$ has $r$ rows, with $r-i+1$ elements in the $i$th row. 
We write $\Li=(a_{i,j})_{\substack {1\leq i\leq r\\ i\leq j\leq r}}$ and draw the pattern aligned to the right:
\begin{equation}
\ytableausetup{mathmode, boxsize=2em}
\begin{ytableau}
a_{1,1} & a_{1,2} &  \cdots  & a_{1,r} \\
\none & a_{2,2} &  \cdots  & a_{2,r} \\
\none & \none &  \ddots  & \vdots \\
\none & \none  &  \none &  a_{r,r} 
\end{ytableau}
\end{equation}

\subsubsection{Type $B_r$ and $C_r$} In this case $\ell\left(\w_0^{B_r}\right)-\ell\left(\w_0^{B_{r-1}}\right)=2r-1.$ These Littelmann patterns have $r$ rows as well, but now the $i$th row has $2r-1$ entries, denoted $a_{i,j}$ for $i\leq j\leq 2r-i.$ We write $\bar{j}=2r-j$ and $\bar{a}_{i,j}=a_{i,\bar{j}},$ and draw the patterns centered as follows: 
\begin{equation}
\ytableausetup{mathmode, boxsize=2em}
\begin{ytableau}
a_{1,1} & a_{1,2} &  \cdots  & a_{1,r}& \cdots & \bar{a}_{1,2} & \bar{a}_{1,1} \\
\none & a_{2,2} &  \cdots  & a_{2,r} & \cdots & \bar{a}_{2,2} & \none \\
\none & \none &  \ddots  & \vdots & \udots  & \none &  \none  \\
\none & \none  &  \none &  a_{r,r} & \none & \none & \none
\end{ytableau}
\end{equation}

\subsubsection{Type $D_r$} In this case $\ell\left(\w_0^{D_r}\right)-\ell\left(\w_0^{D_{r-1}}\right)=2r-2$ for $r\geq 3.$ The $\ell\left(\w_0^{D_r}\right)=r^2-r$ integers from a $BZL$ path are now arranged into a Littelmann pattern with $r-1$ rows. The $i$th row has $2r-2i$ entries, $a_{i,j}$ for $i\leq j\leq 2r-1-i.$ We use notation similar to type $B$ and $C$ and write $\bar{j}=2r-1-j$ for $\bar{a}_{i,j}=a_{i,2r-1-j}.$ 
\begin{equation}
\ytableausetup{mathmode, boxsize=2em}
\begin{ytableau}
a_{1,1} & a_{1,2} &  \cdots  & a_{1,\bar{r}}& a_{1,r}& \cdots & \bar{a}_{1,2} & \bar{a}_{1,1} \\
\none & a_{2,2} &  \cdots  & a_{2,\bar{r}}& a_{2,r} & \cdots & \bar{a}_{2,2} & \none \\
\none & \none &  \ddots  & \vdots & \vdots & \udots  & \none &  \none  \\
\none & \none  &  \none & a_{\bar{r}, \bar{r}}& a_{\bar{r},r} & \none & \none & \none
\end{ytableau}
\end{equation}

\subsection{The bijection with crystal elements}\label{subsect:bijection}

We are ready to give the bijection between crystal elements and Littelmann patterns. 


Recall that the BZL path of a crystal element $b$ consists of $\ell(\w_0)$ segments. Taking the length of these segments produces a tuple $BZL(b)=(n_1,\ldots ,n_{\ell(\w_0)}).$ The entries of the Littelmann pattern $\Li(b)$ corresponding to $b$ are these integers $n_h$ ($1\leq h\leq \ell(\w_0)$). The pattern $\Li(b)$ is filled with elements of $BZL(b)$ row by row proceeding from left to right and from bottom to top. For example, in type $A_r$ we have that $\Li(b)=(a_{i,j})_{\substack {1\leq i\leq r\\ i\leq j\leq r}}$ and:
$$a_{r,r}=n_1,\ a_{r-1,r-1}=n_2,\ a_{r-1,r}=n_3,\ldots ,\ a_{1,1}=n_{\binom{r}{2}+1}, \ldots ,\ a_{1,r}=n_{\binom{r+1}{2}}$$
The shape of the Littelmann patterns above arranges entries in the same column if they correspond to the same edge label. We examine this property in more detail.

\subsection{The weight of a pattern}\label{subsect:wt_of_pattern}

Let $b$ be an element in a crystal element of highest weight $\lambda .$ Let $BZL(b)=(n_1,\ldots ,n_{\ell(\w_0)})$ and $\Li=\Li(b)$ be the Littelmann pattern corresponding to $b$ via the bijection above. Then the weight of $b$ is easily recovered from entries of the pattern. Recall that the $h$th segment of the $BZL$ path follows edges of the crystal labeled with index $k=\w_0(h).$ These edges all correspond to a root operator for the simple root $\alpha_{k}$ i.e. they all have the same label $k$. It follows that 
\begin{equation}\label{eq:wt_diff_BZL}
\lambda -\wt(b)=\sum _{k=1}^r \alpha_k\cdot \sum_{\w_0(h)=k} n_h.
\end{equation}

\begin{figure}[h!]
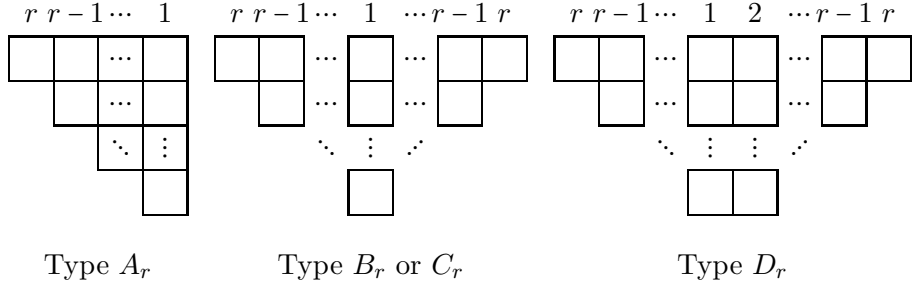

\begin{equation*}
\begin{array}{ccc}
\ytableausetup{mathmode, boxsize=1.5em}
\begin{ytableau}
\none[r] & \none[r-1]  &  \none[\cdots]  & \none[1]\\
 &  &  \cdots  &  \\
\none &  &  \cdots  &  \\
\none & \none &  \ddots  & \vdots \\
\none & \none  &  \none &   
\end{ytableau} & 
\begin{ytableau}
\none[r] & \none[r-1]  &  \none[\cdots]  & \none[1] &  \none[\cdots] & \none[r-1]  & \none[r] \\
 &   &  \none[\cdots]  &    & \none[\cdots] &   &  \\
\none &   &  \none[\cdots]  &  & \none[\cdots] &  & \none \\
\none & \none &  \none[\ddots]  & \none[\vdots]  & \none[\udots]  & \none &  \none  \\
\none & \none  &  \none &   & \none & \none & \none
\end{ytableau} & 
\begin{ytableau}
\none[r] & \none[r-1]  &  \none[\cdots]  & \none[1] & \none[2] & \none[\cdots] & \none[r-1]  & \none[r] \\
 &   &  \none[\cdots]  &  &  & \none[\cdots] &   &  \\
\none &   &  \none[\cdots]  & &  & \none[\cdots] &  & \none \\
\none & \none &  \none[\ddots]  & \none[\vdots] & \none[\vdots] & \none[\udots]  & \none &  \none  \\
\none & \none  &  \none &  &  & \none & \none & \none
\end{ytableau}\\
 & & \\
\text{Type $A_r$} & \text{Type $B_r$ or $C_r$} & \text{Type $D_r$}
\end{array}
\end{equation*}
\caption{Edge colors corresponding to columns of a pattern.}\label{fig:column_index}
\end{figure}

The shape of the patterns has the following property. Entries in a single column of $\Li(b)$ correspond to segments of the BZL path of $b.$ These segments all run along edges of the crystal with the same color $k$ (or $\alpha _k$) $1\leq k\leq r.$ 
Figure \ref{fig:column_index} shows the index of the crystal edges corresponding to each column in the various types. Observe that reading off the index for elements in the top row gives the segment of $\w_0$ that is present in rank $r$ but not in rank $r-1.$

We make this explicit for each of the infinite families. We define the weight $\s(\Li )=(s_1,\ldots ,s_r)$ of a Littelmann pattern $\Li $ ($s_k=s_k(\Li)$) so that:  
\begin{equation}\label{eq:wtdiff_sLi}
\la - \wt (b)=\sum_{k=1}^{r}s_k(\Li(b) )\cdot \alpha_k
\end{equation}

\subsubsection{Type A}
In this case $\Li(b)=(a_{i,j})_{1\leq i\leq j\leq r}.$ A column consists of entries $a_{1,j},\ldots,a_{j,j}.$ If $BZL(b)=(n_1,\ldots,n_{\binom{r+1}{2}})$ then $a_{i,j}=n_{\binom{r-i}{2}+j-i}$ and the corresponding segment of the BZL path of $b$ lies along edges labeled $r-j+1.$ We define $s_k(\Li)$ for $1\leq k\leq r$ by:
\begin{equation}\label{eq:def_patt_wt_comp_A}
s_k(\Li )=\sum_{i=1}^{r}a_{i,r+1-k}
\end{equation}

\subsubsection{Type B or C}
Here we have $\Li(b)=(a_{i,j})_{1\leq i\leq j\leq 2r-i}.$ For any $1<k\leq r$ there are two columns corresponding to the edge index $k,$ the one with $j=r-k+1$ and the one with $\bar{j}=r-k+1.$ We thus define:
\begin{equation}\label{eq:def_patt_wt_comp_BC_not1}
s_k(\Li )=\sum_{i=1}^{r}(a_{i,r+1-k}+\bar{a}_{i,r+1-k})
\end{equation}
The middle column corresponds to the index $1$ and so we define:
\begin{equation}\label{eq:def_patt_wt_comp_BC_1}
s_1(\Li )=\sum_{i=1}^{r}a_{i,r}
\end{equation}
Note that $|\alpha _1|$ is different from $|\alpha_2|=\cdots =|\alpha_r|.$

\subsubsection{Type D}
This case is similar to the previous one. We have $\Li(b)=(a_{i,j})_{1\leq i\leq j\leq 2r-1-i}.$ For $2<k\leq r$ the two columns corresponding to the edge index $k$ are the $j$th where $j=r-k+1$ and $\bar{j}$th, where recall that $\bar{j}=r-k.$ 
The two middle columns corresponds to $\alpha _1$ and $\alpha _2,$ the roots on the ``branched'' end of the Dynkin diagram. Hence we define (cf. \cite{chinta2012littelmann}):
\begin{equation}\label{eq:def_patt_wt_comp_D}
s_k(\Li )=\left\lbrace \begin{array}{ll}
\sum_{i=1}^{r}(a_{i,r+1-k}+\bar{a}_{i,r+1-k}) & \text{if }2<k\leq r\\
\sum_{i=1}^{r}a_{i,r-2+k} & \text{if }k=1,2
\end{array}\right.
\end{equation}


\subsection{Constraints on Littelmann patterns}\label{subsect:Littelmann_inequalities}

The correspondence $b\mapsto \Li (b)$ described in \ref{subsect:bijection} above is a bijection between integral points of the Littelmann cone (see \ref{subsubsect:BZLpaths}) and the set of Littelmann patterns whose entries satisfy a certain set on inequalities, depending on the Cartan type of the underlying root system. To get a set of patterns in bijection with the integral points of a Littelmann polytope for $\lambda$ (equivalently, a crystal of highest weight $\lambda $), we may impose a further set of inequalities on the entries. This second set of inequalities shall depend on the highest weight $\lambda .$ In this section we make these constraints explicit for each of the infinite families of Cartan types. The constructions in section \ref{sect:constructions} phrase the contribution of a crystal element $b$ in terms of whether these inequalities are satisfied by the entries of $\Li(b)$ strictly or with an equality.

\subsubsection{Constraints for the cone}

We give the inequalities describing Littelmann patterns corresponding to integral points of the Littelmann cone $C.$ Let $C^{X_r}$ denote the Littelmann cone in the Cartan type $X_r$. Then we have the following. 

\begin{nthm}\label{thm:Littelman_cone_ineqs}
\cite{littelmann1998cones} Let $b$ correspond to $\Li (b)$ under the bijection described in \ref{subsect:bijection}.
Then $b$ is an integral point of $C^{X_r}$ if and only if the entries of $\Li (b)$ are nonnegative and the following holds. 
\begin{enumerate}
\item[$X_r=A_r:$] \cite[Theorem 5.1]{littelmann1998cones} The rows are weakly decreasing:
\begin{equation}\label{eq:cone_ineq_A}
a_{i,i}\geq a_{i,i+1}\geq \cdots \geq a_{i,r-1}\geq a_{i,r}\text{ for every }1\leq i\leq r
\end{equation}
\item[$X_r=B_r:$] \cite[Theorem 6.1]{littelmann1998cones} For every row we have: 
\begin{equation}\label{eq:cone_ineq_B}
2a_{i,i}\geq 2a_{i,i+1}\geq \cdots \geq 2a_{i,r-1}\geq a_{i,r}\geq 2\bar{a}_{i,r-1} \geq \cdots \geq 2\bar{a}_{i,i}\text{ for every }1\leq i\leq r
\end{equation}
\item[$X_r=C_r:$] \cite[Theorem 6.1]{littelmann1998cones} The rows are weakly decreasing: 
\begin{equation}\label{eq:cone_ineq_C}
a_{i,i}\geq a_{i,i+1}\geq \cdots \geq a_{i,r-1}\geq a_{i,r}\geq \bar{a}_{i,r-1} \geq \cdots \geq \bar{a}_{i,i}\text{ for every }1\leq i\leq r
\end{equation}
\item[$X_r=D_r:$] \cite[Theorem 7.1]{littelmann1998cones} For every row we have:
\begin{equation}\label{eq:cone_ineq_D}
a_{i,i}\geq a_{i,i+1}\geq \cdots \geq a_{i,r-2}\geq a_{i,r-1}, a_{i,r}\geq \bar{a}_{i,r-2} \geq \cdots \geq \bar{a}_{i,i}\text{ for every }1\leq i\leq r-1
\end{equation}
i.e. the rows are weakly decreasing with the exception of the central two elements. There is no restriction on the comparative size of these two elements. 
\end{enumerate}
\end{nthm}

\subsubsection{Constraints for a polytope}

We introduce shorthand to refer to the sums of particular groups of elements of a Littelmann pattern. The notation $s_{i,j}(\Li),$ $\bar{s}_{i,j}(\Li),$ $t_{i,r-1}(\Li),$ $t_{i,r}(\Li)$ used here differs slightly from that of \cite{littelmann1998cones} ($s(a_{i,j}),$ $s(\bar{a}_{i,j}),$ etc.) to emphasize that $s_{i,j}(\Li)$ may be nonzero even if $a_{i,j}$ is not an element of the pattern $\Li .$ When the pattern $\Li $ is clear from context, we write $s_{i,j}$ for $s_{i,j}(\Li).$ 

We define the following shorthand:
\begin{eqnarray}
s_{i,j}(\Li):= &\left\lbrace\begin{array}{ll}
\sum_{k=1}^{i}a_{k,j} & \text{if $\Li $ is type $A$}\\
\sum_{k=1}^{i}(a_{k,j}+\bar{a}_{k,j}) & \text{ if $\Li $ is type $B$ or $C,$ $j\leq r-1$}\\
\sum_{k=1}^{i} a_{k,r} & \text{if $\Li $ is type $B$ and $j=r$}\\
\sum_{k=1}^{i} 2a_{k,r} & \text{if $\Li $ is type $C$ and $j=r$}\\
\sum_{k=1}^{i}(a_{k,j}+\bar{a}_{k,j}) & \text{if $\Li $ is type $D,$ $j\leq r-2$}\\
\sum_{k=1}^{i}(a_{k,r-1}+\bar{a}_{k,r}) & \text{if $\Li $ is type $D,$ $j=r-1$ or $j=r$}
\end{array}\right.\\ \label{eq:def:sij_Li}
\bar{s}_{i,j}(\Li):= &\left\lbrace\begin{array}{ll}
\bar{a}_{i,j}+s_{i-1,j}(\Li) & \text{if $\Li $ is type $B$ and $j\leq r-1$}\\
& \text{or $\Li$ is type $C$ and $j\leq r$} \\
& \text{or $\Li $ is type $D$ and $j\leq r-2$} \\
\bar{a}_{i,j}+2s_{i-1,j}(\Li) & \text{if $\Li $ is type $B$ and $j=r$}\\
s_{i,j}(\Li) & \text{if $\Li $ is type $D,$ $j=r-1$ or $j=r$}
\end{array}\right. \label{eq:def:sbarij_Li}
\end{eqnarray}
Observe that $\bar{s}_{i,j}=s_{i,\bar{j}}$ when both are defined; when only one is, we use this to extend the definition. 
For patterns $\Li $ of type $D$ we shall also need:
\begin{equation}
t_{i,r-1}(\Li ):=\sum_{k=1}^{i} a_{k,r-1} \text{ and } t_{i,r}(\Li ):=\sum_{k=1}^{i} a_{k,r} 
\end{equation}

We are ready to state the inequalities characterizing patterns that correspond to the integral points of a Littelmann polytope, or, equivalently, the points of a crystal of highest weight $\lambda $. Let $\lambda =\sum_{k=1}^{r}m_i\cdot \varpi _i$. The integers $m_k$ appear in the inequalities.

For type $A_r,$ the pattern $\Li$ corresponds to $b\in \Cr_{\la }$ if the following inequalities are satisfied \cite[Corollary 4]{littelmann1998cones}: 
\begin{equation}\label{eq:polytope_ineq_A}
a_{i,j}\leq m_{r-j+1}+s_{i,j-1}(\Li )-2s_{i-1,j}(\Li )+s_{i-1,j+1}(\Li );\text{ for $1\leq i\leq j\leq r$}
\end{equation}

For type $B_r$ and type $C_r$ the inequalities are as follows \cite[Corollary 6.]{littelmann1998cones}: 
\begin{eqnarray}\label{eq:polytope_ineq_BC}
\bar{a}_{i,j}\leq & m_{r-j+1}+\bar{s}_{i,j-1}(\Li )-2s_{i-1,j}(\Li )+s_{i-1,j+1}(\Li );\text{ for $1\leq i\leq j\leq r-1$}\\
a_{i,j}\leq & m_{r-j+1}+\bar{s}_{i,j-1}(\Li )-2\bar{s}_{i,j}(\Li )+s_{i,j+1}(\Li );\text{ for $1\leq i\leq j\leq r-1$}\\
a_{i,r}\leq & m_{1}+d\bar{s}_{i,r-1}(\Li )-d\bar{s}_{i-1,r}(\Li )
\end{eqnarray}
where $d=2$ in type $B$ and $d=2$ in type $C.$

Finally for type $D_r$ the inequalities are as follows \cite[Corollary 8.]{littelmann1998cones}:
\begin{eqnarray}\label{eq:polytope_ineq_D}
\bar{a}_{i,j}& \leq & m_{r-j+1}+\bar{s}_{i,j-1}(\Li )-2s_{i-1,j}(\Li )+s_{i-1,j+1}(\Li );\text{$1\leq i\leq j\leq r-2$}\\
a_{i,j}& \leq & m_{r-j+1}+s_{i,j+1}(\Li )-2\bar{s}_{i,j}(\Li )+\bar{s}_{i,j-1}(\Li );\text{$1\leq i\leq j\leq r-2$}\\
a_{i,r-1}& \leq & m_{2}+\bar{s}_{i,r-2}(\Li )-2t_{i-1,r-1}(\Li )\\
a_{i,r}& \leq & m_{1}+\bar{s}_{i,r-2}(\Li )-2t_{i-1,r}(\Li )
\end{eqnarray}

Let $BZL(\lambda )$ denote the set of Littelmann patterns $\Li $ that are in bijection with elements of the highest weight crystal $\Cr_{\lambda }.$ 

\section{The Constructions}\label{sect:constructions}

We are ready to give the constructions of $p$ parts of WMDS and metaplectic Whittaker functions, i.e. the objects from section \ref{subsect:MDS_WF}. To emphasize the combinatorial nature of the constructions in this section we restrict our attention to constructing a polynomial $P_{\lambda }$. The meaning of this polynomial in each of the types was given in \ref{subsect:review_lit} and \ref{subsect:MDS_WF}.

We shall give a polynomial $P$ for any Cartan type in the infinite families $A_r,$ $B_r,$ $C_r$ and $D_r.$ 

The constructions are analogous in different types. Before giving the type by type constructions in \ref{subsect:typebytypeconstr} we begin by summarizing the common elements. In this section and afterwards, we shall identify a crystal element with the corresponding Littelmann pattern. 

\subsection{The contribution of a pattern}\label{subsect:shape_of_contribution}

In all cases $P=P_{\lambda }$ is a sum over a crystal of highest weight $\lambda $. Using the bijection $b\mapsto \Li (b)$ above, we write $P$ as a sum over $BZL(\la ).$ We shall have a sum:
\begin{equation}\label{eq:P_sum_expression}
P=\sum_{\Li \in BZL(\la )} G(\Li)\cdot \x^{\wt(\Li)}
\end{equation}
where recall that $\CC[\Lambda ]$ was identified with a polynomial ring $\CC[\x].$ 
Here $wt(\Li )$ is essentially the weight of the pattern given in \ref{subsect:wt_of_pattern}. 

The coefficient $G(\Li )$ shall be given as a product:
\begin{equation}\label{eq:G_is_product}
G(\Li )=\prod _{i,j} g_{i,j}(\Li )
\end{equation}
This product is over elements of the Littelmann pattern, and the factor $g_{i,j}(\Li )$ depends only on the {\em{decoration}} of the element $a_{i,j}$ in $\Li .$ For each of the infinite families, we decorate the elements of the pattern $\Li $ according to a {\em{circling}} and a {\em{boxing}} rule. An entry $a_{i,j}$ may be circled, boxed, neither, or both. Before giving the rules for decorating elements of $\Li $ in \ref{subsect:decoration_rules}, we preview how the decorations affect the contribution of $\Li .$

In the case of the constructions in type $A,$ $B$ and $C,$ the factor $g_{i,j}(\Li )$ depends only on the decoration of $a_{i,j}$ and the integer $a_{i,j}$ itself.\footnote{In type $D$ the picture is more complex. An analogous construction gives a result slightly different from the one expected from the $p$-part. The phenomenon of {\em{(symmetric) multiple leaners}} accounts for this discrepancy; see \ref{subsubsect:TypeD_constr} for details.} 
In particular (in type $A$ and $C$)\footnote{Note that in type $B$ we have the coefficient $1$ instead of $q^{a_{i,j}}.$ This discrepancy by a factor of $q,$ also present when comparing  \eqref{eq:coeffs_typeB} with \eqref{eq:coeffs_typeA} or \eqref{eq:coeffs_typeC} can be eliminated with a change of variables in the polynomial $P_{\la}(\x).$}:
\begin{equation}
\text{If $a_{i,j}$ is circled, then }g_{i,j}(\Li ):=\left\lbrace\begin{array}{ll}
q^{a_{i,j}} & \text{if $a_{i,j}$ is not boxed}\\
0 & \text{if $a_{i,j}$ is boxed}
\end{array}\right.
\end{equation}
If $a_{i,j}$ is {\em{not}} circled, then the value of $g_{i,j}(\Li )$ is a Gauss sum (see \ref{subsect:Gauss_sums}).

\subsection{Circling and boxing rules}\label{subsect:decoration_rules}

Recall that in \ref{subsubsect:BZLpaths} we identified elements of the highest weight crystal $\Cr _{\la}$ with integral points of the Littelmann polytope. In section \ref{sect:Littelmann_patterns} we gave a set of constraints that a pattern $\Li (b)\in BZL(\la )$ satisfies if it corresponds to an element $b\in \Cr _{\la}.$ The constraints came in the form of inequalities \eqref{eq:cone_ineq_A}-\eqref{eq:cone_ineq_D} (these guarantee that $b$ is in the Littelmann cone) and \eqref{eq:polytope_ineq_A} - \eqref{eq:polytope_ineq_D} (these are specific to the polytope and depend on $\la $). The decoration of an entry $a_{i,j}$ of $\Li $ depends on whether the inequalities involving $a_{i,j}$ are satisfied with an equality. 

\subsubsection{Circling rule}\label{subsubsect:circling}

The inequalities \eqref{eq:cone_ineq_A}-\eqref{eq:cone_ineq_D} involve a single row of the Littelmann pattern. An element $a_{i,j}$ appears in one of these, and that has a lower bound for $a_{i,j}.$ The lower bound is of the form:
\begin{equation}\label{eq:ineq_aij_circ}
a_{i,j}\geq \left\lbrace\begin{array}{ll}
\max(a_{i,j+1}, a_{i,j+2}) & \text{if $\Li $ is of type $D$ and $j=r-2$}\\
a_{i,j+2} & \text{if $\Li $ is of type $D$ and $j=r-1$}\\
\frac{1}{2}\cdot a_{i,j+1} & \text{if $\Li $ is of type $B$ and $j=r-1$}\\
2\cdot a_{i,j+1} & \text{if $\Li $ is of type $B$ and $j=r$}\\
a_{i,j+1} & \text{otherwise}
\end{array}\right.
\end{equation}

Then $a_{i,j}$ is {\em{circled}} if this lower bound \eqref{eq:ineq_aij_circ} holds with an equality. 

\subsubsection{Boxing rule}\label{subsubsect:boxing}

The inequalities \eqref{eq:polytope_ineq_A} -\eqref{eq:polytope_ineq_D} have the property that every $a_{i,j}$ entry of a Littelmann pattern appears on the left-hand side of exactly one of them. Let $a_{i,j}$ be {\em{boxed}} if this inequality holds with an equality. 

\subsubsection{Interpretation of decorations}\label{subsubsect:interpret_decorations}

Observe that an entry is circled or boxed when an inequality defining the Littelmann polytope is satisfied with an equality. This means that the element of the polytope corresponding to $\Li $ is on one of the hyperplanes defining the polytope. 

\subsection{Constructions type by type}\label{subsect:typebytypeconstr}

We are ready to finish describing the constructions in type $A_r,$ $B_r,$ $C_r$ and $D_r.$ 

\subsubsection{Type A}\label{subsubsect:TypeA_constr}

We recall the definition of a $p$-part from \cite{bbf-wmdbook}. First recall that the weight $\s(\Li )$ of a pattern was defined in \eqref{eq:def_patt_wt_comp_A} as the sum of entries in columns. In this case we are interested in the $p$-part \eqref{eq:MDS-ppart}; we set $P=P_{\la}(\x)$ equal to this $p$-part with $\la =\sum_{i=1}^{r} (l_i+1)\varpi _i$ and we assume that $\la $ is a strongly dominant weight ($m_i=l_i\geq 1$ for every $1\leq i\leq r$). 

Then $P$ is given by \eqref{eq:P_sum_expression} where $\wt(\Li )=\s(\Li )$ and $G(\Li )$ is a product as in \eqref{eq:G_is_product} where the factors of the coefficient are given as follows \cite[Chapter 1]{bbf-wmdbook}: 
\begin{equation}\label{eq:coeffs_typeA}
g_{i,j}(\Li )=\left\lbrace \begin{array}{ll}
0 &\text{if $a_{i,j}$ is circled and boxed}\\
q^{a_{i,j}}&\text{if $a_{i,j}$ is circled but not boxed}\\
g(a_{i,j}) &\text{if $a_{i,j}$ is boxed but not circled}\\
h(a_{i,j}) &\text{if $a_{i,j}$ is neither circled nor boxed}
\end{array}\right.
\end{equation}
Here $g=g_1$ and $h=h_1$ are the Gauss sums from \ref{subsect:Gauss_sums} and $q$ is the order of a residue field.

\begin{nrem}
This is the construction of the $p$ part $H_{\Gamma}$ from \cite{bbf-wmdbook}. 
\end{nrem}

\subsubsection{Type B}\label{subsubsect:TypeB_constr}
Recall from \ref{subsect:review_lit} and \ref{subsubsect:Whittaker_fns} that \cite{brubaker2012metaplectic} gives a conjectural formula for a Whittaker function in type $B.$ In this case we have that $P_{\la }(\x)$ is the value $\W(\x,\la)$ of a Whittaker function on a torus element. 

In this case let $\wt(\Li ):=\sum_{k=1}^{r}s_k\alpha_k-\lambda $ where $\s(\Li)=(s_1,\ldots ,s_r)$ as defined in \eqref{eq:def_patt_wt_comp_BC_not1} and \eqref{eq:def_patt_wt_comp_BC_1}. Then \cite[Conjecture 2]{brubaker2012metaplectic} states that $P$ is given by \eqref{eq:P_sum_expression} and \eqref{eq:G_is_product} where the factors are given as follows: 
\begin{equation}\label{eq:coeffs_typeB}
g_{i,j}(\Li )=\left\lbrace \begin{array}{ll}
0 &\text{if $a_{i,j}$ is circled and boxed}\\
1 &\text{if $a_{i,j}$ is circled but not boxed}\\
q^{-a_{i,j}}g_t(a_{i,j}) &\text{if $a_{i,j}$ is boxed but not circled}\\
q^{-a_{i,j}}h_t(a_{i,j}) &\text{if $a_{i,j}$ is neither circled nor boxed}
\end{array}\right.
\end{equation}
where $g_t$ and $h_t$ are as in \ref{subsect:Gauss_sums} and the subscript $t$ is $1$ if $j=r$ and $t=2$ otherwise. 

\subsubsection{Type C}\label{subsubsect:TypeC_constr}

Once again we have $P_{\la }(\x)$ be the $p$-part of a Multiple Dirichlet series from \cite{beineke2012crystal} or \cite{friedberg2015eisenstein}. Let us once again write $\wt(\Li )=\s(\Li ).$ To specify $P $ by \eqref{eq:P_sum_expression} and \eqref{eq:G_is_product} we must again specify the factors $g_{i,j}(\Li ):$
\begin{equation}\label{eq:coeffs_typeC}
g_{i,j}(\Li )=\left\lbrace \begin{array}{ll}
q^{a_{i,j}} &\text{if $a_{i,j}$ is circled but not boxed}\\
g_t(a_{i,j}) &\text{if $a_{i,j}$ is boxed but not circled}\\
h_1(a_{i,j}) &\text{if $a_{i,j}$ is neither circled nor boxed and $n|a_{i,j}$}\\
0 &\text{otherwise}
\end{array}\right.
\end{equation}
where again $t=1$ if $j\neq r$ and $t=2$ is $j=r.$ 
(We note that \cite[(31)]{beineke2012crystal} contained a typo that was fixed by \cite[(34)]{friedberg2015eisenstein}: note that by \eqref{subsect:Gauss_sums} if $a_{i,j}$ is neither circled nor boxed, then $g_{i,j}(\Li )=0$ unless $n|a_{i,j}.$)

\subsubsection{Type D}\label{subsubsect:TypeD_constr}
Finally, we recall \cite[Conjecture 1]{chinta2012littelmann}, a conjectural expression for the $p$-part $P_{\la }(\x)$  of a Multiple Dirichlet series. 
In this case the construction for $P(=P_{\la })$ is slightly different. Once again it is a sum \eqref{eq:P_sum_expression} over contributions from Littelmann patterns $\Li \in BZL(\la )$ and the weight of a pattern is $\wt(\Li )=\s(\Li ).$ The contribution $G(\Li )$ of a pattern is again a product. However, in this case the factor $g_{i,j}$ is dependent on more than the decoration of $a_{i,j}.$ 

The definition of $G(\Li )$ in \cite{chinta2012littelmann} is written as a product over connected components of the {\em{decorated graph}} $\Gamma (\Li )$ of the pattern $\Li .$ To give the conjectural construction of the $p$-part, we introduce some terminology. 

The vertices of $\Gamma (\Li )$ are the entries of $\Li .$ Two entries belong to the same connected component in $\Gamma (\Li )$ if they are comparable in the inequalities \eqref{eq:cone_ineq_D} and they are equal. By the rightmost element of a component $C$ me mean the entry in $C$ that is positioned rightmost in the Littelmann pattern. A connected component $C$ is called a {\em{multiple leaner}} (m.l.) if it consists of entries $a_{i,j_1}=a_{i,j_1+1}=\cdots =a_{i,j_2}$ where $j_1\leq r-2,$ $r+1\geq j_2$ and $a_{i,j_1-1}>a_{i,j_1},$ $a_{i,j_2}>a_{i,j_2+1}.$ By the legs of $C$ we mean the entries $a_{i,j_1}=\cdots =a_{i,r-2}$ and $a_{i,r+1}=\cdots =a_{i,j_2};$ the entry on the endpoint of the shorter leg of $\Gamma $ is $a_{i,j_1}$ or $a_{i,j_2}.$ The component $C$ is called a {\em{symmetric multiple leaner}} (s.m.l.) if in addition  $j_2=\bar{j}_1;$ in this case we define its length to be $l(C)=r-j_1,$ half the number of its vertices. 

We may write $G(\Li )$ as in \eqref{eq:G_is_product} but to define $g_{i,j}(\Li )$ we write \cite[5.5]{chinta2012littelmann} 
\begin{equation}\label{eq:typeD_prod_over_components}
\sigma (C)=\prod _{a_{i,j}\in C} g_{i,j}(\Li )
\end{equation}
and give $\sigma (C)$ in terms of standard contributions of the entries it contains. Let \cite[5.5]{chinta2012littelmann} 
\begin{equation}\label{eq:sigma_entry}
\sigma (y)=\left\lbrace \begin{array}{ll}
0 & \text{if the entry $a$ is circled and boxed}\\
h_1(a)\cdot q^{-a} & \text{if the entry $a$ is not boxed and not circled}\\
g_1(a)\cdot q^{-a} & \text{if the entry $a$ is boxed and not circled}
\end{array}\right.
\end{equation}
We then define $\sigma(C)$ to be as follows:
\begin{itemize}
\item $\sigma(C)=0$ if any $a_{i,j}\in C$ is both circled and boxed;
\item $\sigma(C)=\sigma (a)$ if $C$ is not a m.l. and $a$ is its rightmost element, or $C$ is a m.l. that is not symmetric and $a$ is the endpoint of its shorter leg; 
\item $\sigma(C)=\sigma (a)(1-q^{-l(C)})$ if $C$ is a s.m.l., $a\neq 0$ is its rightmost element, and $a$ is unboxed;
\item $\sigma(C)=\sigma (a_{i,j})\sigma(a_{i,j-1})q^{1-l(C)} $ if $C$ is a s.m.l., $a_{i,j}$ is its rightmost element, and $a_{i,j}\neq 0$ is boxed;
\item $\sigma(C)=1$ if $C$ is a s.m.l. with zero entries.
\end{itemize}
Note that $\sigma (C)$ is a product over $g_{i,j}(\Li ),$ but now $g_{i,j}(\Li)$ depends not only on the decoration of $a_{i,j}\in \Li ,$ but also on the position of $a_{i,j}$ within a connected component of $\Gamma (\Li ),$ and whether that component is a (symmetric) multiple leaner or not. 

\section{Branching}\label{sect:branching_properties}

In the previous section we described constructions of polynomials $P_{\la }(\x)$ that are of interest from a number theoretic perspective as explained in \ref{subsect:MDS_WF}. We also mentioned in \ref{subsect:motivation_branching} that elucidating the relationship of the constructions with the branching properties of highest weight crystals can be the key to understanding some of their properties. In this section, we take a closer look at how the branching properties of crystals manifests in these constructions. 

In all of the examples above the polynomial $P$ was associated to a crystal $\Cr _{\la }$ corresponding to a root system $\Phi $ or rank $r,$ with Cartan type in one of the infinite families $A_r,$ $B_r,$ $C_r$ or $D_r.$ When the edges of $\Cr_{\la }$ labeled by $\alpha _r$ are omitted, the remaining graph is a disjoint union of rank $r-1$ crystals of the same Cartan type, but rank $r-1:$ 
\begin{equation}\label{eq:crystal_branching}
\Cr _{\la }=\bigsqcup_{\mu } \Cr _{\mu }
\end{equation}
We now wish explain how $P_{\lambda}$ can be written in terms of the polynomials $P_{\mu}$ associated to those crystals. 

The construction of $P_{\la }$ is given (cf. \eqref{eq:P_sum_expression}) as a sum over Littelmann patterns $\Li \in BZL(\la ),$ each contributing a term $G(\Li )\x^{\wt(\Li )}.$ We examine how the sets $BZL(\mu )$ corresponding to $\mu $ in \eqref{eq:crystal_branching} can be recovered from $BZL(\la )$ by giving a pattern $\Li '\in BZL(\mu )$ for any pattern $\Li \in BZL(\la )$ in \ref{subsect:BZLbranching}. We indicate a method of computing the weights $\mu $ that appear in \eqref{eq:crystal_branching} in \ref{subsect:wts_appearing}. 
Then in \ref{subsect:wts_appearing} and \ref{subsect:contrib_decomp} we explain how for such pairs $\Li $ and $\Li '$ the contributions $G(\Li )$ and $\wt(\Li )$ can be written in terms of $G(\Li '),$ $\wt (\Li ')$ and $\mu .$ For the remainder of the discussion let us fix a dominant weight $\la =\sum_{k=1}^{r}m_k\varpi _k$ and the corresponding crystal $\Cr _{\la }.$

\subsection{Patterns with fixed top row}\label{subsect:BZLbranching}

Recall the bijection between elements of the crystal $\Cr _{\la }$ and Littelmann patterns given in \ref{subsect:bijection}. For an element $b\in \Cr _{\la },$ the entries of the pattern $\Li (b)=\Li _{\la }(b)$ are the lengths of the segments in the BZL path of $b.$ Here $BZL(b)$ corresponds to the choice of a particular long word $\w_0.$ 

Let the element $b\in \Cr_\la $ belong to $\Cr _{\mu }$ in the decomposition \eqref{eq:crystal_branching}. We shall sometimes write $b'$ when we mean $b$ as an element of the abstract crystal $\Cr _{\mu };$ write $\Li '=\Li (b')\in BZL(\mu ).$  

As remarked in \ref{subsect:choice_longword} the choices made in \eqref{eq:choice_longword} all have the property that the long word $\w_0^{X_r}$ chosen in rank $r$ starts with the long word $\w_0^{X_{r-1}}$ chosen in rank $r-1.$ Together with the shape of the patterns (see \ref{subsect:shape_of_patterns}) this means that $\Li (b')$ is the same as $\Li (b)$ without its first row. This argument proves the following. 

\begin{nlem}\label{lem:BZLbranching}
Let $b_\mu \in \Cr _{\mu} \subset \Cr _{\la }$ be the highest element within a crystal in the decomposition \eqref{eq:crystal_branching}. Let $\Li _{\mu }(b)$ and $\Li _{\la }(b)$ denote the Littelmann patterns corresponding to any $b\in \Cr _{\mu} \subset \Cr _{\la }$ as an element of $\Cr _{\mu}$ and $\Cr _{\la },$ respectively. Then for any $b\in \Cr _{\mu} \subset \Cr _{\la }$ the top row of $\Li _{\la }(b)$ is the same as the top row of $\Li _{\la }(b_{\mu }),$ and $\Li _{\mu }(b)$ can be recovered from $\Li _{\la }(b)$ by deleting the top row. 
\end{nlem}

\subsection{The weights in the decomposition}\label{subsect:wts_appearing}
In light of Lemma \ref{lem:BZLbranching} and the inequalities on the top row of $\Li \in BZL(\la )$
we can describe the highest weights $\mu $ appearing in the decomposition \eqref{eq:crystal_branching} by computing the weight of the highest element $b_\mu \in \Cr _{\mu}.$ Recall that by \eqref{eq:wtdiff_sLi} $\la -\wt(b_{\mu })=\la -\mu $ can be expressed in terms of $\s(\Li _{\la }(b_{\mu })):$ 
\begin{equation}\label{eq:lambdamudiff_from_sLibmu}
\la -\mu =\sum_{k=1}^{r} \s_k (\Li _{\la }(b_{\mu }))\cdot \alpha _k
\end{equation}
and furthermore that we have:
\begin{equation}\label{eq:s_branching}
\s(\Li _{\la }(b))=\s(\Li _{\la }(b_{\mu }))+\s(\Li _{\mu }(b)) \text{ for any $b\in \Cr_{\mu} \subset \Cr _{\la }$}
\end{equation}
Note that since the entries under the first row of $\Li _{\la }(b_{\mu })$ are all zero, the right hand side of \eqref{eq:lambdamudiff_from_sLibmu} can be written entirely in terms of the entries in the first row of $\Li _{\la }(b_{\mu }).$ The inequalities restricting the first row of a $\Li \in BZL(\la )$ involve no entries from any other row (cf. \ref{subsect:Littelmann_inequalities}). It follows that given a highest weight $\la ,$ we can recover the set of weights $\mu $ that appear in the decomposition \eqref{eq:crystal_branching}. (This involves expressing the simple roots $\alpha _k$ in terms of the fundamental weights, and carefully examining the restrictions on entries of the first row of a pattern $\Li \in BZL(\la )$.) We omit further discussion of this here and refer the reader to \cite[(2.4)]{bbf-wmdbook} for an example of a similar statement in Cartan type $A.$

\subsection{Branching and contributions}\label{subsect:contrib_decomp}

Let $b$ be an element $b\in \Cr _{\mu} \subset \Cr _{\la }.$ Let $\Li =\Li _{\la}(b)$ and $\Li ' =\Li _{\mu}(b).$ We wish to write $P_{\la }(\x)$ as a sum 
\begin{equation}\label{eq:P_branched}
P_{\la }(\x) =\sum _{\mu } p(\mu )\cdot P_{\mu}(\x )
\end{equation}
where the weights $\mu $ are the ones of the decomposition \eqref{eq:crystal_branching}, and $p(\mu )$ is a monomial. 

{Before we explain why this decomposition is possible for the polynomials $P_{\la }(\x)$ we remark on the terminology of ``branching.'' Recall that the polynomial $P_{\la }(\x)$ can be thought of as a deformation of a highest-weight character. Equation \eqref{eq:P_branched} has a clear analogue for highest weight characters. Let us write $V_{\la}$ and $V_{\mu}$ for the irreducible representations of highest weight $\lambda $ and $\mu $ (in rank $r$ and $r-1$) respectively. If $P_{\la }(\x)$ were the character of $V_{\la }$ then the coefficient in the monomial $p(\mu )$ would match the multiplicity of the $V_\mu $ in the restriction of $V_{\la }$ to a subalgebra of corank one determined by the first $r-1$ simple roots.}

The contribution of a pattern $\Li $ is of the form $G(\Li )\x^{\wt(\Li )}$ as in \eqref{eq:P_sum_expression}. The term $\x^{\wt(\Li )}$ depends only on $\s(\Li ).$ It follows from \eqref{eq:s_branching} that the first $r-1$ component of $\s(\Li )$ is the tuple $\s(\Li '),$ while $\s_k(\Li )=\s_k(\Li (b_{\mu }))$ depends only on $\mu .$ 

We turn next to the coefficient $G(\Li ).$ Recall from section \ref{subsect:shape_of_contribution} and in particular \eqref{eq:G_is_product} that $G(\Li )$ is a product of factors $g_{i,j}(\Li ).$ The factor $g_{i,j}(\Li )$ essentially depends on the decoration of an entry $a_{i,j}$ in $\Li ,$ and the decorations in turn depend on whether the inequalities imposed on the entry $a_{i,j}$ by $\Li $ being an element of $BZL(\la )$ are satisfied strictly or with an equality. It is immediate that the factors $g_{1,j}(\Li )$ corresponding to entries of the first row of $\Li $ depend only on $\mu .$ Closer examination of the inequalities imposed on the lower rows and the decorations show that in fact $g_{i,j}(\Li )$ can be written as a product of $g_{i-1,j-1}(\Li ')$ and a factor that depends only on $\mu ,$ and not the element $b\in \Cr _{\mu }.$ 

Thus we may conclude that one may in fact decompose $P_\la (\x)$ as in \eqref{eq:P_branched}. For a precise statement of this flavor in Cartan type $A,$ see \cite[Proposition 16.]{puskas2016whittaker}.

\bibliographystyle{amsalpha_no_mr_initials.sty}
\bibliography{crystals_in_NT}

\end{document}